\documentclass{amsart}
\usepackage{amssymb,latexsym,amsmath,amsthm}
\usepackage[mathscr]{eucal}
\newtheorem{theorem}{Theorem}[section]
\newtheorem{lemma}[theorem]{Lemma}

\numberwithin{equation}{section}

\renewcommand{\l}{\lambda}
\newcommand{\g}{{\rm g}}
\newcommand{\RR}{\ensuremath{\mathbb{R}}}
\newcommand{\CC}{\ensuremath{\mathbb{C}}}
\newcommand{\ZZ}{\ensuremath{\mathbb{Z}}}
\newcommand{\sph}{\ensuremath{\mathbb{S}}}
\newcommand{\prtl}{\ensuremath{\partial}}

\newcommand{\hf}{\ensuremath{\frac{1}{2}}}
\newcommand{\supp}{\ensuremath{{\rm supp}}}

\title
[Strichartz estimates for Schr\"{o}dinger operators] {On
Strichartz estimates for Schr\"{o}dinger operators in compact
manifolds with boundary}

\thanks{The authors were supported by the National Science Foundation,
Grants DMS-0140499, DMS-0354668, DMS-0555162, and DMS-0354???.}

\author{Matthew D. Blair}
\address{Department of Mathematics, Johns Hopkins University,
Baltimore, MD 21218}
\email{mblair@math.jhu.edu}

\author{Hart F. Smith}
\address{Department of Mathematics, University of Washington,
Seattle, WA 98195}
\email{hart@math.washington.edu}

\author{Christopher D. Sogge}
\address{Department of Mathematics, Johns Hopkins University,
Baltimore, MD 21218}
\email{sogge@jhu.edu}

\begin{document}

\maketitle

\section{Introduction}
Let $(M,\g)$ be a Riemannian manifold of dimension $n \geq 2$.
Strichartz estimates are a family of dispersive estimates on
solutions $u(t,x): [-T,T] \times M \to \CC$ to the Schr\"{o}dinger
equation
\begin{align}\label{E:schrodman}
D_t u + \Delta_\g u &= 0, & & u(0,x) = f(x)
\end{align}
where $\Delta_\g$ denotes the Laplace-Beltrami operator on $(M,\g)$
and $D_t= i^{-1}\partial_t$. In their most general form, local
Strichartz estimates state that
\begin{equation}\label{E:strichartz}
\|u\|_{L^p([-T,T]; L^q(M))} \leq C\|f\|_{H^s(M)}
\end{equation}
where $H^s$ denotes the $L^2$ Sobolev space over $M$, and $2 \leq
p,q \leq \infty$ satisfies
\begin{align}\label{E:stzpair}
\frac{2}{p} + \frac{n}{q} = \frac{n}{2}, & & (n,p,q)\neq
(2,2,\infty).
\end{align}
Such estimates are well established in flat Euclidean space, where
$M=\RR^n$ and $\g_{ij}=\delta_{ij}$.
In that case $s=0$, and one can take $T=\infty$; see for example
Strichartz~\cite{strich77}, Ginibre and Velo~\cite{ginvelo85},
Keel and Tao~\cite{keeltao98}, and references therein. Estimates
for the standard flat 2-torus were shown by Bourgain~\cite{Bourgain93}
to hold for any $s>0$.

There is also considerable interest in developing these estimates
for non-flat geometries, and also for compact domains.
In the case where $M$ is compact and $\prtl M =
\varnothing$, Burq, G\'{e}rard, and Tzvetkov \cite{burq1}
established~\eqref{E:strichartz} with $s=\frac{1}{p}$. Hence there
is a loss of derivatives in their estimate when compared to
the case of flat geometries.

A simple investigation
of the Schr\"{o}dinger evolution on spherical harmonics where
$M=\sph^n$ shows that some loss of derivatives must occur.
For instance, with $n=2$, by taking the initial data to be a
highest weight spherical harmonic on $\sph^2$ one concludes that
the best possible local $L^2_x\to L^4_{t,x}$ bounds would involve a loss
of $1/8$ derivatives.  This sharp estimate and related ones for
Zoll surfaces were obtained in \cite{burq1}.  It is not known, however, whether
the weaker estimates involving a loss of $\frac 1p$ derivatives
in \cite{burq1} for general compact manifolds without
boundary can be improved.

In the case where $\prtl M \neq \varnothing$, one also considers
Dirichlet or Neumann boundary conditions in addition
to~\eqref{E:schrodman}
\begin{align*}
u(t,x)|_{\prtl M} =0 \text{   (Dirichlet)}, & & N_x \cdot \nabla
u(t,x)|_{\prtl M} =0 \text{   (Neumann)}
\end{align*}
where $N_x$ denotes the unit normal vector field to $\prtl M$.
Here one expects a further loss of derivatives in the estimates. The
Rayleigh whispering gallery modes over the unit disk in $\RR^2$
provide examples of Dirichlet eigenfunctions which accumulate
their energy near the boundary, contributing to high $L^p$ norms.
Applying the Schr\"{o}dinger evolution to these eigenfunctions
show that $s\ge\frac{1}{6}$ is necessary for the
Strichartz estimate with $p=q=4$.
Recently, Anton~\cite{Anton} showed that the
estimates~\eqref{E:strichartz} hold on general manifolds
with boundary provided $s >\frac{3}{2p}$.
In addition, the arguments of~\cite{Anton} work equally well for
a manifold $M$ without boundary equipped with a
metric $\g$ of Lipschitz regularity.

In this work, we improve on the current results for compact
$(M,\g)$ where either $\prtl M \neq \varnothing$,
or $\prtl M = \varnothing$
and $\g$ is Lipschitz, by showing that Strichartz estimates hold with a
loss of fewer derivatives.

\begin{theorem}\label{T:c1strich}
Let $(M,\g)$ be either a smooth compact Riemannian manifold with boundary,
or a manifold without boundary equipped with
a Lipschitz metric $\g$.  Then the following Strichartz
estimate holds for any Strichartz pair~\eqref{E:stzpair}
\begin{equation}\label{E:c1strich}
\|e^{it\Delta_\g} f \|_{L^p([-T,T]; L^q(M))} \lesssim
\|f\|_{H^{\frac{4}{3p}}(M)}
\end{equation}
\end{theorem}

In the case where $(M,\g)$ is a  boundaryless manifold with $\g \in
C^\infty$, the estimate of
Burq-G\'{e}rard-Tzvetkov~\eqref{E:strichartz} with
$s=\frac{1}{p}$, while not known to be sharp, is a natural
result by the following heuristic argument. For a general
compact manifold, there are no conjugate points for the geodesic flow
at distance less than the injectivity radius of the manifold.
Given a solution to the Schr\"{o}dinger equation whose
frequencies are concentrated at $\lambda$, energy propagates
at speed $\approx\lambda$. Hence, a frequency $\lambda$ solution should
possess good dispersive properties at least until time
$T_\lambda \approx \frac{1}{\lambda}$. We thus expect to be able
to prove a Strichartz estimate with no loss of derivatives for
such a solution over a time interval of size roughly
$\frac{1}{\lambda}$. By considering a sum over such intervals we
should obtain a Strichartz estimate over a time interval of unit
size, only with a constant appearing on the right hand size
which is a constant multiple of $\lambda^\frac{1}{p}$.  This
corresponds to $s=\frac{1}{p}$ in the estimate, and
Littlewood-Paley theory yields the estimate for arbitrary solutions.

In the case where $\prtl M \neq \varnothing$, the boundary conditions
affect the flow of energy near the boundary.  A
key strategy involves reflecting the metric and the solution across
the boundary, to obtain a Schr\"odinger equation on
a manifold without boundary, but with a metric that has Lipschitz
singularities along $\prtl M$. Hence matters reduce to
considering the Schr\"{o}dinger evolution for Lipschitz metrics.
In this case, when establishing estimates for solutions at frequency
$\lambda$, one can replace the rough metric by a regularized
metric which has conjugate points at distance roughly
$\lambda^{-\frac{1}{3}}$ apart. Therefore, the solutions should
possess good dispersive properties over a time interval of size
roughly $\lambda^{-\frac{4}{3}}$. This now yields a Strichartz
estimate over a time interval of unit size with a loss of
$\frac{4}{3p}$ derivatives.  Hence, for manifolds with boundary
\eqref{E:c1strich} appears to be the natural analog of the
aforementioned estimates of \cite{burq1} for the general
boundaryless case.

Our proof of Theorem~\ref{T:c1strich} follows the above heuristics.
In section 2
the solution is localized spatially and a coordinate chart is used
to work on $\RR^n$; a
Littlewood-Paley decomposition then reduce matters to establishing
Strichartz estimates for components of the solutions dyadically
localized in frequency. As alluded to above, we then seek to prove
Strichartz estimates with no loss of derivatives over time
intervals of size $\lambda^{-\frac{4}{3}}$ for components of the
solution localized at frequency $\lambda$.  This involves
regularizing the metric by truncating its frequency to a scale
dependent on $\lambda$. Rescaling the solution then reduces the
problem to establishing Strichartz estimates for metrics with 2 bounded
derivatives over small time intervals whose size also depends on
the frequency.  Section 3 uses a phase space transform to
construct a parametrix for such Schr\"{o}dinger operators, and
section 4 concludes the paper by showing that the parametrix
yields the desired estimates.

\subsection*{Notation}
In what follows $d$ will denote the gradient operator which maps
scalar functions to vector fields and vector fields to matrix
functions in the natural way. The expression $X \lesssim Y$ means
that $X \leq CY$ for some $C$ depending only on $n$ and on the
Lipschitz norm of the metric.

\section{Reductions}

We will establish Theorem~\ref{T:c1strich} more generally for
operators on $M$ which take the following form in local coordinates
\begin{equation}\label{Pform}
\bigl(Pf\bigr)(x)=
\rho(x)^{-1}\sum_{i,j=1}^n
\partial_i\Bigl(\rho(x)\,\g^{ij}(x)\,\partial_jf(x)\Bigr)
\,.
\end{equation}
Such an operator is self-adjoint in the measure
$d\mu=\rho(x)\,dx$. Neumann conditions and the boundary normal are
defined with respect to the metric $\g_{ij}$.

We start by reducing the case of a
manifold $M$ with boundary and $P$ smooth, to the case of a compact manifold
$M$ without boundary, with $P$ having coefficients of Lipschitz regularity.
For this, let $\tilde M$ denote the double of $M$, identified along $\prtl M$.
We define a differentiable structure on $\tilde M$ near $\prtl M$
using geodesic normal coordinates in $\g_{ij}$, so $x_n>0$ and $x_n<0$ define
the two copies of $M$. In these coordinates, $\g^{ni}=0$ for $i\ne n$,
hence $P$ contains no cross terms between $\partial_n$ and $\partial_i$.
The operator $\tilde P$ with
coefficients $\g^{ij}(x',|x_n|)$ and $\rho(x',|x_n|)$ is thus symmetric
under $x_n\rightarrow -x_n$, and extends the lift of $P$ to $\tilde M$
across $\prtl M$ to one with Lipschitz coefficients.
Eigenspaces for $\tilde P$ decompose into
symmetric and antisymmetric functions; these
correspond to extensions of
eigenfunctions for $P$ satisfying
Dirichlet (resp. Neumann) conditions, and each eigenfunction
is of regularity $C^{1,1}$
across the boundary. The Schr\"odinger flow for $\tilde P$ is thus
easily seen to extend
that for $P$, and Strichartz estimates for $P$ follow by establishing
such estimates for $\tilde P$ on $\tilde M$.

We assume henceforth that $M$ is a compact manifold with smooth differentiable
structure, on which an operator $P$ of the form \eqref{Pform} is given,
with coefficients of Lipschitz regularity.
Define $L^q$-Sobolev spaces on $M$ using the spectral resolution
of $P$,
$$
\|f\|_{W^{s,q}(M)}=\|\langle D_P\rangle^s f\|_{L^q(M)}\,,
\qquad \langle D_P\rangle=\bigl(1-P\bigr)^{\frac 12}\,.
$$
When $q=2$ we denote $W^{s,q}$ by $H^s$.
By elliptic regularity (e.g. \cite[Theorem 8.10]{GT}
for $q=2$, and  \cite[Theorem 9.11]{GT} or
\cite[\S 2.2]{Tay} for other $q$) the spaces $W^{s,q}$ for $1<q<\infty$
coincide with the Sobolev spaces defined
using local coordinates, provided $0\le s\le 2$.

Suppose that $u(t,x) = (e^{itP} f) (x)$.
Then we need to establish
\begin{equation*}
\|u\|_{L^p([-T,T];W^{s,q}(M))}
\lesssim\|f\|_{H^1(M)}\,,\qquad s=1-\tfrac 4{3p}\,.
\end{equation*}
Let $u=\sum_{k=0}^\infty u_k$ denote
a Littlewood-Paley partition of $u$, where $u_k$ for $k\ge 1$ is
spectrally localized to
$\langle D_P\rangle\approx 2^k$.
Then, for $p,q\ge 2$,
\[
\|\langle D_P\rangle^s u\|_{L^p_tL^q_x}\approx
\|(\langle D_P\rangle^s u)_k\|_{L^p_t L^q_x\ell^2_k}\le
\|\langle D_P\rangle^s u_k\|_{\ell^2_kL^p_t L^q_x}\,,
\]
and
\[
\|f\|_{H^1(M)}\approx \|f_k\|_{\ell^2_k H^1(M)}\,,
\]
hence it suffices to show, uniformly over $k$, that
\begin{equation*}
\|u_k\|_{L^p([-T,T];W^{s,q}(M))}
\lesssim\|f_k\|_{H^1(M)}\,,\qquad s=1-\tfrac 4{3p}\,.
\end{equation*}
By taking a finite partition of
unity, it suffices to prove that
\begin{equation}\label{E:halfspace}
\| \psi u_k \|_{L^p([-T,T];W^{s,q}(\RR^n))} \lesssim
\|u_k \|_{L^\infty([-T,T];H^1(M))}
\end{equation}
for each smooth cutoff $\psi $ supported in a suitably
chosen coordinate chart. We will choose
coordinate charts such that the image contains the unit ball,
and
$$
\|\g^{ij}-\delta_{ij}\|_{Lip(B_1(0))}\le c_0\,,
\qquad
\|\rho-1\|_{Lip(B_1(0))}\le c_0\,,
$$
for $c_0$ to be taken suitably small. (This may require multiplying
$\rho$ by a harmless constant). We take $\psi$ supported in the
unit ball, and assume $\g^{ij}$ and $\rho$ are extended so that
the above holds globally on $\RR^n$.

Let $\{\beta_j(D) \}_{j\geq 0}$ be a Littlewood-Paley partition of
unity on $\RR^n$, and
$v_j=\beta_j(D)\psi u_k$.
We will prove that, for each $j$,
\begin{equation}\label{vjest}
\|\langle D\rangle^s v_j\|_{L^p_tL^q_x}\lesssim
\|v_j\|_{L^\infty_t H^1_x}+
2^{-\frac j3}\|(D_t+P)v_j\|_{L^\infty_t L^2_x}\,,
\end{equation}
with all norms taken over $[-T,T]\times\RR^n$, and
$\langle D\rangle=(1-\Delta)^{\frac 12}$.

This will imply \eqref{E:halfspace}, provided we dominate the sum over
$j$ of the right hand side of \eqref{vjest}
by $\|u_k\|_{L^\infty([-T,T];H^1(M))}$, which we now do.

For a Lipschitz function $a$, $\;[\beta_j(D),a]:H^{s-1}\rightarrow H^s$
for $s=0,1$. Hence $[P,\beta_j(D)\psi]:H^1\rightarrow L^2$,
and it follows that
$$
\|(D_t+P)v_j\|_{L^\infty_tL^2_x}\lesssim
\|u_k\|_{L^\infty_t H^1(M)}\,,
$$ hence the second term on the right of \eqref{vjest}
is bounded by a geometric series.
For the first term, note that
\begin{align*}
\|v_j\|_{L^\infty_tH^1_x}&\lesssim
\min\bigl(2^j\|v_j\|_{L^\infty_tL^2_x},2^{-j}\|v_j\|_{L^\infty_tH^2_x}\bigr)\\
&\lesssim
\min\bigl(2^j\|u_k\|_{L^\infty_tL^2(M)},2^{-j}\|u_k\|_{L^\infty_tH^2(M)}\bigr)
\,.
\end{align*}
The sum over $j$ is dominated by
$$
2^k\|u_k\|_{L^\infty_tL^2(M)}+2^{-k}\|u_k\|_{L^\infty_tH^2(M)}
\lesssim
\|u_k\|_{L^\infty_tH^1(M)}\,.
$$

Setting $\l=2^j$, and denoting $v_j$ by $v_\l$,
the estimate \eqref{vjest} is equivalent to
\begin{multline*}
\|v_\l\|_{L^p([-T,T];L^q(\RR^n))}\\
\lesssim
\l^{\frac 4{3p}}\|v_\l\|_{L^\infty([-T,T];L^2(\RR^n))}+
\l^{\frac 4{3p}-\frac 43}\|(D_t+P)v_\l\|_{L^\infty([-T,T];L^2(\RR^n))}
\,.
\end{multline*}
This, in turn, follows by showing that for any interval $I_\l$ of
length $\l^{-\frac 43}$, we have
\begin{equation}\label{vlest}
\|v_\l\|_{L^p(I_\l;L^q(\RR^n))}\\
\lesssim
\|v_\l\|_{L^\infty(I_\l;L^2(\RR^n))}+
\|(D_t+P)v_\l\|_{L^1(I_\l;L^2(\RR^n))}
\,.
\end{equation}

We now regularize the coefficients of $P$ by setting
$$
\g^{ij}_\l=S_{\l^{-2/3}}(\g^{ij})\,,\qquad
\rho_\l=S_{\l^{-2/3}}(\rho)\,,
$$
where $S_{\l^{2/3}}$ denotes a truncation of a function to
frequencies less than $\l^{\frac 23}$, and let $P_\l$
denote the operator with coefficients $\g^{ij}_\l$ and $\rho_\l$.
Since $|\g^{ij}_\l-\g^{ij}|\lesssim \l^{-\frac 23}$,
and similarly for $\rho$, it follows that
$$
\|(P-P_\l)v_\l\|_{L^1(I_\l;L^2(\RR^n))}\lesssim
\|v_\l\|_{L^\infty(I_\l;L^2(\RR^n))}\,,
$$
and we may thus replace $P$ by $P_\l$ on the right hand side
of \eqref{vlest} without changing the estimate.

Finally, we rescale the problem. Let $\mu=\l^{\frac 23}$, and
define
\begin{align*}
u_\mu(t,x) &=v_{\lambda}(\l^{-\frac 23}t,
\l^{-\frac 13} x)\\
Q_\mu (x,D) &= P_\lambda(\lambda^{-\frac{1}{3}}x,D).
\end{align*}
The function $u_\mu(t,\cdot)$ is localized to frequencies of size
$\mu$, and the coefficients of $Q_\mu$ are localized to
frequencies of size less than $\mu^{\frac 12}$. This implies the
following estimates on the coefficients of $Q_\mu$
$$
|\prtl^\alpha_x \g^{ij}_\lambda (\l^{-\frac 13}x)|+
|\prtl^\alpha_x \rho_\lambda (\l^{-\frac 13}x)|\leq
C_\alpha \,\mu^{\frac 12 \max(0, |\alpha| -2)} .
$$
The interval $I_\l$ scales to an interval of length $\mu^{-1}$.
We have thus reduced the proof of Theorem~\ref{T:c1strich} to the
following.

\begin{theorem}\label{T:c2strich}
Suppose that $u_\mu(t,x)$ is localized to frequencies $|\xi| \in
[\frac{1}{4} \mu, 4 \mu]$ and solves
\begin{equation}\label{E:schrodeqn}
\Bigl(\,D_t + \sum_{1 \leq i,j \leq n} a^{ij}_\mu(x)\,\prtl_{x_i}
\prtl_{x_j} +\sum_{1 \leq i \leq n} b^i_\mu(x)\, \prtl_{x_i}\,\Bigr) u_\mu
= F_\mu
\end{equation}
Assume also that the metric satisfies
\[
\|a^{ij}_\mu-\delta_{ij}\|_{C^2} \ll 1\,,\qquad
\|b^i_\mu\|_{C^1} \lesssim 1\,,
\]
\[
\supp\bigl(\,\widehat{a^{ij}_\mu}\,\bigr)\,,
\supp\bigl(\widehat{b_\mu^i}\,\bigr) \subset
B_{\mu^{1/2}}(0)\,.
\]
Then the following estimate holds
$$
\|u_\mu\|_{L^p([0,\mu^{-1}];L^q(\RR^n)}) \lesssim
\|u_\mu\|_{L^\infty([0,\mu^{-1}];L^2(\RR^n))} +
\|F_\mu\|_{L^1([0,\mu^{-1}];L^2(\RR^n))}\,.
$$
\end{theorem}

\section{The Parametrix}

We will establish Theorem~\ref{T:c2strich} using a short-time wave packet
parametrix for the equation \eqref{E:schrodeqn}. Wave packet parametrices
have been used to establish Strichartz estimates for Schr\"odinger
equations in the work of Staffilani-Tataru \cite{ST} and Koch-Tataru \cite{KT};
see Tataru \cite{T} for an overview of the methods.
The result we need, in fact, is
included as a special case in Theorem 2.5 of \cite{KT}. The proof of
the short time estimate Theorem~\ref{T:c2strich} is comparatively simple,
though, and therefore we include a self-contained proof here
for the reader's benefit.

In this section, then, we use a wave packet transform to construct a
parametrix for the operator
$$
D_t + A(x,D) + B(x,D)=D_t + \sum_{1 \leq i,j \leq n}
a^{ij}_\mu(x)\prtl_{x_i} \prtl_{x_j} +\sum_{1 \leq i \leq n}
b^i_\mu(x) \prtl_{x_i}
$$
in~\eqref{E:schrodeqn} that will yield the Strichartz estimates.
For convenience, we suppress the $\mu$ from both the operator and
the solution.  Let $g$ be a radial Schwartz function over $\RR^n$
such that $\supp(\widehat{g}) \subset B_1(0)$ and $\|g\|_{L^2} =
(2\pi)^{-\frac{n}{2}}$.  For $\mu \geq 1$, we define the operator
$T_\mu: \mathcal{S}'(\RR^n) \to C^\infty (\RR^{2n})$ by
$$
T_\mu f(x, \xi) = \mu^{\frac{n}{4}} \int e^{-i\langle \xi, z-x
\rangle}g(\mu^\hf(z-x))f(z)dz.
$$
$T_\mu$ enjoys the property that its adjoint as a map from
$L^2_{x,\xi}(\RR^{2n}) \to L^2_{z}(\RR^n)$ also serves as a left
inverse for $T_\mu$, that is, $T_\mu^* T_\mu=I$.  This implies
that $T_\mu$ is an isometry
$$
\|T_\mu f\|_{L^2_{x,\xi}(\RR^{2n})} = \|f\|_{L^2_z(\RR^n)}.
$$

We conjugate $A(x,D)$ by $T_\mu$ and take a suitable
approximation to the resulting operator.  Specifically, define the
following differential operator over $(x,\xi)$
$$
\widetilde{A} = - i d_\xi a(x,\xi) \cdot d_x + i d_x
a(x,\xi) \cdot d_\xi  + a(x,\xi) - \xi \cdot d_\xi a(x,\xi)\,.
$$
By a standard argument from wave packet methods
(see for example \cite{T1} or \cite {T} where $g$ is Gaussian,
or Lemmas 3.1-3.3 in~\cite{SmC2} for $g$ as above)
we have that if $\widetilde{\beta}_\mu$ is a
Littlewood-Paley cutoff truncating to frequencies $|\xi| \approx
\mu$ then
\begin{equation}\label{E:conjug}
\|T_\mu A(\cdot,D) \widetilde{\beta}_\mu(D) - \widetilde{A} T_\mu
\widetilde{\beta}_\mu(D)\|_{L^2_z \to L^2_{x,\xi}} \lesssim \mu\,.
\end{equation}

This yields that, if $\tilde{u} (t,x,\xi) = (T_\mu
u_\mu(t,\cdot))(x,\xi)$, then $\tilde{u}$ solves the equation
$$
\bigl(\prtl_t + d_\xi a(x,\xi) \cdot d_x - d_x a(x,\xi) \cdot d_\xi + i
a(x,\xi) - i\xi\cdot d_\xi a(x,\xi)\bigr) \tilde u(t,x,\xi) =
\tilde{G}(t,x,\xi)\,,
$$
where $\tilde{G}$ satisfies
$$
\int_0^{\mu^{-1}} \|\tilde{G}(t,x,\xi)\|_{L^2_{x,\xi}} dt \lesssim
\|u_\mu\|_{L^\infty([0,\mu^{-1}];L^2)} +
\|F_\mu\|_{L^1([0,\mu^{-1}];L^2)}.
$$

Given an integral curve $\gamma(r) \in \RR^{2n}_{x,\xi}$ of the
vector field
$$
\prtl_t + d_\xi a(x,\xi) \cdot d_x - d_x a(x,\xi) \cdot d_\xi
$$
with $\gamma(t) = (x,\xi)$,  we denote $\chi_{s,t}(x,\xi) =
(x_{s,t},\xi_{s,t})= \gamma(s)$.  Now define
\begin{align*}
\sigma (x,\xi) = a(x,\xi) - \xi\cdot d_\xi a(x,\xi) \,,& &
\psi(t,x,\xi) = \int_0^t \sigma(\chi_{r,t}(x,\xi)) \;dr.
\end{align*}
This allows us to write
\begin{equation*} \tilde{u} (t,x,\xi) =
e^{-i\psi(t,x,\xi)}u_{0}(\chi_{0,t}(x,\xi)) + \int_0^t
e^{-i\psi(t-r,x,\xi)}\,
\tilde{G}(r,\chi_{r,t}(x,\xi))\,dr\,,
\end{equation*}
which expresses $\tilde u$ as an integrable superposition over $r$
of functions invariant under the flow of $\tilde A$, truncated
to $t>r$.

Since $u(t,x) = T_\mu^* \tilde{u}(t,x,\xi)$ it thus suffices to
obtain estimates
$$
\|\widetilde{\beta}_\mu (D)W_t f\|_{L^p_t L^q_x} \lesssim
\|f\|_{L^2_{x,\xi}}
$$
where $W_t$ acts on functions $f(x,\xi)$ by the formula
$$
(W_t f)(y) = T_{\mu}^*\bigl(e^{-i\psi(t,x,\xi)}f(\chi_{0,t}(x,\xi))\bigr)(y)\,.
$$
By a standard duality argument and an application of the
endpoint estimates of Keel-Tao~\cite{keeltao98} this results from
establishing
\begin{equation}\label{E:dispersive}
\|\widetilde{\beta}_\mu W_t W_s^* \widetilde{\beta}_\mu \|_{L^1
\to L^\infty} \lesssim |t-s|^{-\frac{n}{2}}
\end{equation}
\begin{equation}\label{E:energy}
\|\widetilde{\beta}_\mu W_t W_s^* \widetilde{\beta}_\mu \|_{L^2
\to L^2} \lesssim 1
\end{equation}
The inequality \eqref{E:energy} follows from the fact that $T_\mu$ is an
isometry and $\chi_{0,t}(x,\xi)$ is a symplectomorphism, hence
preserves the measure $dx\,d\xi$.  The
inequality \eqref{E:dispersive} is the focus of the next section.

\section{The dispersive estimate}
In this section, we establish the inequality~\eqref{E:dispersive}.
We write the kernel $K(t,y,s,x)$ of $W_t W_s^*$ as
$$
\mu^{\frac{n}{2}} \int e^{-i\langle \zeta, x-z \rangle - i
\int_s^t \sigma(\chi_{r,t}(z,\zeta))\;dr +i \langle \zeta_{t,s},
y-z_{t,s} \rangle} g(\mu^\hf(y-z_{t,s}))g(\mu^\hf(x-z))\,dz\, d\zeta\,.
$$
Recall that $\supp(\hat{g}) \subset B_1(0)$. Since we are
concerned with $W_t W_s^* \widetilde{\beta}_\mu$, we can insert a
cutoff $S_\mu(\zeta)$ into the integrand which is supported in a
set $|\zeta| \approx \mu$.  Note that the
Hamiltonian vector field is independent of time, and hence
$\chi_{t,s}= \chi_{t-s,0}$.  We drop the zero and abbreviate
the latter transformation as $\chi_{t-s}(z,\zeta) =
(z_{t-s},\zeta_{t-s})$. It then suffices to consider $s=0$,
and we write the kernel $K(t,x,0,y)$ as
$$
\mu^{\frac{n}{2}} \int e^{-i\langle \zeta, x-z \rangle - i
\psi(t,z,\zeta) +i \langle \zeta_t,
y-z_t \rangle} g(\mu^\hf(y-z_t))g(\mu^\hf(x-z))S_\mu (\xi)
\,dz\, d\zeta\,.
$$
We need to establish uniform bounds over $x$ and $y$,
$|K(t,x,0,y)| \lesssim t^{-\frac{n}{2}}$.
A straightforward estimate shows that
$$
|K(t,x,0,y)| \lesssim \mu^n
$$
meaning that the dispersive estimate holds for $t \leq \mu^{-2}$.
We thus assume $t\ge \mu^{-2}$ for the remainder of the
section.  Lastly we suppose that, in addition, $t \leq
\varepsilon \mu^{-1}$ with $\varepsilon$ chosen sufficiently
small and independent of $\mu$.

We first need derivative estimates on the transformation
$\chi_t(z,\zeta)$.
\begin{lemma}
Consider the solutions $(z_t(z,\zeta),\zeta_t(z,\zeta))$ to
Hamilton's equations
\begin{equation}\label{E:hameqn}
\prtl_t z_t = d_\xi a_\mu(z,\zeta)\,,\qquad
\prtl_t \zeta_t = -d_x a_\mu(z,\zeta)\,,\qquad
(z_0\,,\zeta_0)=(z\,,\zeta)\,.
\end{equation}
We then have the following estimates on the first partial
derivatives of $(z_t,\zeta_t)$ when $|\zeta| \in [\frac{1}{4}\mu,
4\mu] $ and $|t| \leq \mu^{-1}$
\begin{equation}\label{E:1derivest}
\begin{matrix}
|d_z z_t -I| \lesssim \mu t &\qquad\qquad& |d_\zeta z_t| \lesssim t\\ \\
|d_z \zeta_t| \lesssim \mu^2  t & &  |d_\zeta \zeta_t -I| \lesssim
\mu t
\end{matrix}
\end{equation}
\medskip
\begin{equation}\label{E:1derivest2}
\bigl|\,d_\zeta z_t - \int_0^t \bigl(d^2_\zeta a_\mu\bigr)(\chi_s(z,\zeta))\,ds\,\bigr|
\lesssim \mu t^2
\end{equation}
The higher partial derivatives satisfy, for $j+k\ge 2$,
\begin{align}\label{E:kderivest}
\mu\,|d^j_z d^k_\zeta z_t | + | d^j_z d^k_\zeta \zeta_t | \lesssim
\mu^{2-k}\,t\,\langle \mu^{\frac{3}{2}}t \rangle^{j+k-1}.
\end{align}
\end{lemma}
\begin{proof}If $|\zeta|\approx 1$, then we can write the Hamilton
equations as:
$$
(z_t,\zeta_t)=(z,\zeta)+\int_0^t v(z_s,\zeta_s)\,ds\,,
$$
where the vector field $v$ satisfies
$$
|d_{z,\zeta}^k v|\lesssim \mu^{\frac 12(k-1)}\,,\qquad k\ge 1\,.
$$
Differentiating the equation and using induction yields
the bound,
$$
|d_{z,\zeta}^k(z_t,\zeta_t)-d_{z,\zeta}^k(z,\zeta)|\lesssim
t\,\langle\mu^{\frac 12} t\rangle^{k-1}\,,\qquad |t|\le 1\,.
$$
Estimates \eqref{E:1derivest} and \eqref{E:kderivest} now follow by
the rescaling property
$$
(z_t(z,\zeta),\zeta_t(z,\zeta)) = \bigl(z_{\mu t}(z,\mu^{-1}\zeta)\,,
\mu\,\zeta_{\mu t}(z,\mu^{-1}\zeta)\bigr)\,.
$$
Estimate \eqref{E:1derivest2} follows by differentiating Hamilton's
equations as above and applying the bounds \eqref{E:1derivest}.
\end{proof}

We take a partition of unity $\{\phi_m\}_{m\in\ZZ^n}$ over $\RR^n$
with
$\phi_m(\zeta) = \phi(t^{\hf}(\zeta-t^{-\hf}m))$ for some $\phi$
smooth and compactly supported.  We then write
$$
K(t,y,0,x) = \sum_{m\in \ZZ^n} K_m(t,y,x)
$$
where $K_m(t,y,x)$ is defined by
$$
\mu^{\frac{n}{2}}\! \int e^{-i\langle \zeta, x-z \rangle
-i\psi(t,z,\zeta) +i \langle \zeta_{t},
y-z_{t} \rangle} g(\mu^\hf(y-z_{t}))\,g(\mu^\hf(x-z))\,\phi_m (\zeta)
S_\mu (\zeta)\, dz\, d\zeta
$$
The key estimate is that, for $\xi_m = t^{-\hf} m$,
\begin{equation}\label{E:kest}
|K_m(t,x,y)| \lesssim
t^{-\frac{n}{2}}\bigl(1+t^{-\hf}|y-x_t(x,\xi_m)|\,\bigr)^{-N}.
\end{equation}
Estimate \eqref{E:1derivest2}, and the fact that
$$
\|d_\xi^2 a(x,\xi)-2I\|=2\|a^{ij}-\delta^{ij}\|\ll 1
$$
yields for $l,m \in \ZZ^n$ and $t \leq \varepsilon\mu^{-1}$
$$
|x_t(x,\xi_m)-x_t(x,\xi_l)| \approx t\,|\xi_m-\xi_l| = t^\hf |m-l|\,.
$$
This now yields
$$
\sum_{m\in \ZZ^n} |K_m(t,x,y)| \lesssim t^{-\frac{n}{2}}
\sum_{m\in \ZZ^n}(1+|m|)^{-N}.
$$
Since the sum on the right converges for $N$ large this
establishes the dispersive estimate.

To prove \eqref{E:kest}, we start by noting that
$$
\prtl_{\zeta_i}\left(\int_0^{t} a(z_{r},\zeta_{r}) -
\zeta_r\cdot(d_\zeta a)(z_{r},\zeta_{r})\,dr \right) + \zeta_t
\cdot \prtl_{\zeta_i} z_t=0.
$$
The expression vanishes
at $t=0$ since $d_\zeta z_0=0$, and Hamilton's equations
show that the derivative of the expression with respect to
$t$ vanishes.

As in Theorem 5.4 of Smith-Sogge~\cite{SmSoBdry}, we now proceed
by defining the differential operator
$$
L = \frac{1+it^{-1}(x-z-d_\zeta \zeta_t \cdot (y - z_t))\cdot
d_\zeta }{1+t^{-1}|x-z-d_\zeta \zeta_t \cdot (y - z_t)|^2}.
$$
By the observation above, $L$ preserves the phase function in
the definition of $K_m$.
The estimates \eqref{E:1derivest} and \eqref{E:kderivest} show that,
if $p$ is any one of the functions $\phi_m(\zeta)$, $t^{-\frac 12}z_t$,
$\mu^{\frac 12}z_t$, $S_\mu(\zeta)$, $\mu^{-\frac 12}t^{-\frac 12}\zeta_t$,
then for $\mu^{-2}\le t\le \mu^{-1}$,
$$
|(t^{-\frac 12}\partial_\zeta)^k p|\lesssim 1\,.
$$
Integration by parts now yields the following upper bound on
$K_m(t,x,y)$
\begin{multline*}
\mu^{\frac{n}{2}} \int_{\RR^n \times \supp(\phi_m)}
(1+t^{-1}|(x-z)-d_\zeta \zeta_t \cdot (y -
z_t)|^2)^{-N}\\\times(1+\mu^{\hf}|x-z|)^{-N}(1+\mu^\hf|y-z_t|)^{N}
\,dz\, d\zeta
\end{multline*}

We conclude by showing that
\begin{equation}\label{E:xzapprox}
t^{-\hf} |(x-z)-d_\zeta \zeta_t \cdot (x_t - z_t)| \lesssim
1+\mu|x-z|^2
\end{equation}
where $x_t$ denotes $x_t(x,\xi_m)$.  This implies that the
integrand is dominated by
$$
(1+t^{-1}|d_\zeta \zeta_t \cdot
(y-x_t)|^2)^{-N}(1+\mu^{\hf}|x-z|)^{-N}\,.
$$
Since $|d_\zeta \zeta_t -I| \lesssim \varepsilon$, this
establishes the estimate \eqref{E:kest}, since the $z$ decay
and compact $\zeta$ support imply that the
integral is essentially over a region in phase space of volume roughly
$t^{-\frac{n}{2}}\mu^{-\frac{n}{2}}$.

To establish~\eqref{E:xzapprox}, we employ a Taylor
expansion and~\eqref{E:kderivest} to obtain
\begin{multline*}
t^{-\hf}|x_t - z_t-(d_z z_t)(x-z)- (d_\zeta z_t)(\xi_m-\zeta)|\\
\lesssim
t^{\frac 12}\langle\mu^{\frac 32}t\rangle
\bigl(\,\mu|x-z|^2+|x-z|\,|\xi_m-\zeta|+\mu^{-1}|\xi_m-\zeta|^2\bigr)
\lesssim 1 +\mu |x-z|^2
\end{multline*}
where the last inequality uses the fact that $\mu^{-2} \leq t\le \mu^{-1}$
and $|\,\xi_m-\zeta| \lesssim t^{-\hf}$.  In addition, by \eqref{E:1derivest}
$$
t^{-\hf}|(d_\zeta z_t)(\xi_m-\zeta)| \lesssim t (t^{-\hf})^2 =1.
$$
Since $\chi_t(z,\zeta)$ is a symplectomorphism, we have
$$
\prtl_{\zeta_i}\zeta_t \cdot \prtl_{z_j} z_t - \prtl_{\zeta_i }
z_t \cdot \prtl_{z_j} \zeta_t=\delta_{ij}
$$
where $\cdot$ pairs the $z_t$ and $\zeta_t$ indices.  Lastly,
by \eqref{E:1derivest},
$$
t^{-\hf}|d_\zeta z_t | \, |d_z \zeta_t| \, |x-z|
\lesssim \mu^2t^{\frac{3}{2}}|x-z| \leq \mu^{\hf}|x-z|\,.
$$
These facts now combine to yield the estimate~\eqref{E:xzapprox}.\qed

\end{document}